\documentclass[
a4paper 
,11pt 
,leqno 
]{scrartcl} 

\usepackage{my-preamble-article} 
\usepackage{my-macros} 

\title{Uniqueness of the pendent drop of infinite length}

\usepackage{authblk} 
\author{Emmanuel \textsc{Risler}\footnote{\protect\url{http://maths.insa-lyon.fr/~erisler/}}}
\affil{Université de Lyon, CNRS, INSA de Lyon\\ Institut Camille Jordan (CNRS UMR 5208)\\ 21 avenue Jean Capelle, 69621 Villeurbanne, France}

\begin{document}
\maketitle
\thispagestyle{empty}
\begin{abstract}
We prove the uniqueness of the infinite length axisymmetric solution to the capillary equation. We observe that capillary equation can be viewed, at large depth, as a perturbation of an integrable two-dimensional differential system. Uniqueness is then proved by an elementary perturbation argument. 
\end{abstract}
\section{Introduction}
An orientable hypersurface is said to satisfy a prescribed mean curvature equation if the mean curvature at each point of the hypersurface is prescribed by a scalar field. When this scalar field is proportional to the ``vertical'' coordinate, this equation is commonly referred as the \emph{capillary equation}, since it is satisfied by a static interface between two liquid phases of different densities, where the capillary force (proportional to the mean curvature) has to balance the difference of pressures (which, up to a change of the origin, is proportional to this vertical coordinate). 

\begin{figure}[!htbp]
	\centering
    \includegraphics[width=0.3\textwidth]{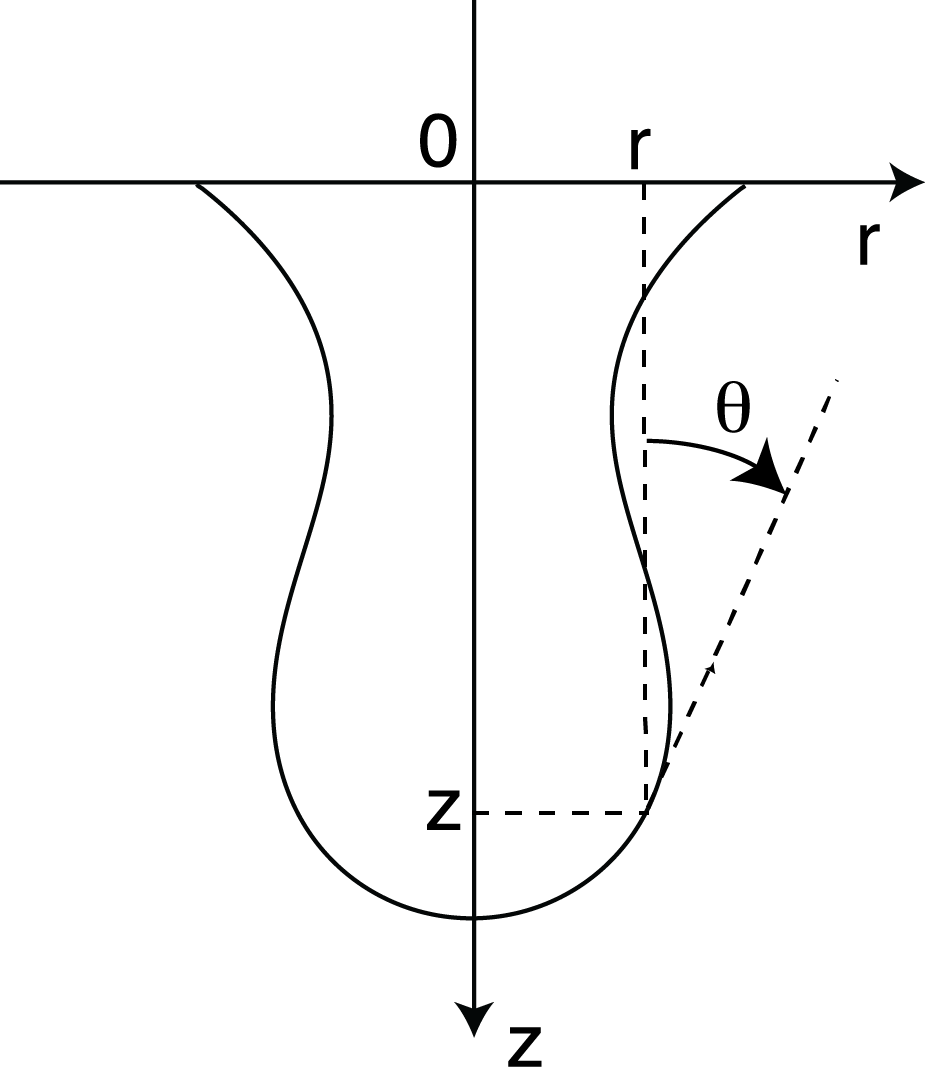}
    \caption{Definition of angle $\theta$.}
    \label{fig_def_theta}
\end{figure}

Let $n$ be an integer not smaller than $2$, and let $(x_1,\dots,x_n,z)$ denote the coordinates in $\rr^{n+1}$ (the ``vertical'' coordinate will be $z$). A hypersurface of $\rr^{n+1}$ is axisymmetric with respect to the $z$-axis if, at every point of the hypersurface, the distance $r=\sqrt{x_1^2+\dots+x_n^2}$ to the $z$-axis is a function of $z$. Let us introduce the angle $\theta$ as displayed on figure~\ref{fig_def_theta} (namely, satisfying $dr/dz = \tan\theta$). Then, at every point of this axisymmetric hypersurface, the capillary equation reads:
\begin{equation}
\label{axisym_cap_equ}
(n-1)\frac{\cos\theta}{r}+\cos\theta\,\frac{d\theta}{dz}=\kappa z\,,
\end{equation}
where $\kappa$ is a fixed real quantity. Let us assume without loss of generality that $\kappa$ is positive (the vertical coordinate $z$ increases ``downwards''). When $n$ equals $2$ equation~(\ref{axisym_cap_equ}) governs the shape of an axisymmetric liquid pendent drop in $\rr^3$. 

Up to a change of scale (in both $r$ and $z$), we can choose whatever positive value for $\kappa$. Let us make the convenient choice $\kappa=n-1$; indeed, with this choice, equation~(\ref{axisym_cap_equ}) can be rewritten as the following two-dimensional differential system:
\begin{equation}
\label{syst_rz}
\left\{
\begin{aligned}
dr/dz &= -\tan\theta \\
d\theta/dz &= (n-1)\Bigl(\frac{z}{\cos\theta}-\frac{1}{r}\Bigr)
\end{aligned}
\right.
\qquad\mbox{where}\quad
r>0
\quad\mbox{and}\quad
\theta\in\Bigl(-\frac{\pi}{2},\frac{\pi}{2}\Bigr)\,.
\end{equation}

Drops of finite length correspond to solutions of equation~(\ref{axisym_cap_equ}) --- or equivalently~(\ref{syst_rz}) --- that do not exist for $z$ larger than a certain value $z_{\max}$ (the height at the tip of the drop) where $r(z)$ vanishes. In this paper we are interested in solutions that are defined up to $+\infty$ in $z$, namely drops of infinite length. The main result is the following. 
\begin{theorem}
\label{main_theorem}
There exists a unique solution $z\mapsto (r,\theta)$ of system~(\ref{syst_rz}) that is  defined on $(0,+\infty)$; it satisfies:
\[
r(z)\sim 1/z
\quad\mbox{and}\quad
\theta(z)\rightarrow0
\quad\mbox{when}\quad
z\rightarrow+\infty\,.
\]
\end{theorem} 
Historically, the problem goes back to P.-S.~Laplace \cite{laplace1805}, T.~Young \cite{young1805}, and Lord Kelvin \cite{kelvin1891}. From a modern mathematical point of view, the main contributions are due to P.~Concus and R.~Finn. In \cite{concusfinn1974,concusfinn1979} (see also Finn's book \cite{finn1986}) they proved the existence of solutions corresponding to drops of every prescribed finite length, and studied their shapes. In \cite{concusfinn1975a}, they proved the existence of a solution under the form of a function $z$ of the radius $r$, for which $z(r)$ approaches $+\infty$ when $r$ approaches $0$, thus corresponding to a drop of infinite length, and computed its asymptotic expansion when $r$ approaches $0$. In \cite{concusfinn1975b} they achieved a first step towards its uniqueness, by proving uniqueness among all singular solutions asymptotically sufficiently close to that asymptotic expansion. Further steps towards uniqueness were achieved by M.-F.~Bidaut-Véron (\cite{bidautveron1986,bidautveron1996}), who proved uniqueness under a weaker explicit criterion on the closeness to the asymptotic expansion. Finally, the full proof of the uniqueness was achieved by R.~Nickolov (\cite{nickolov2002}) who proved that the Bidaut-Véron uniqueness criterion is actually satisfied by every singular solution. 

In \cite{rierarisler2002}, C.~Riera and the author proposed to view and study system~(\ref{syst_rz}) as a dynamical system. This point of view, which surprisingly enough does not seem to have been adopted before, turned out to be fruitful: classical dynamical systems techniques led us to describe the dynamics of all solutions, and to recover in a simple way many results on the shape of pendent drops of finite length, together with the existence of (at least) one singular solution corresponding to a drop of infinite length. The aim of the present paper is to provide, by the same approach, a short, elementary, and self-contained proof of its uniqueness. 
\section{Rescaling and statement of main result}
\begin{figure}[!htbp]
	\centering
    \includegraphics[width=\textwidth]{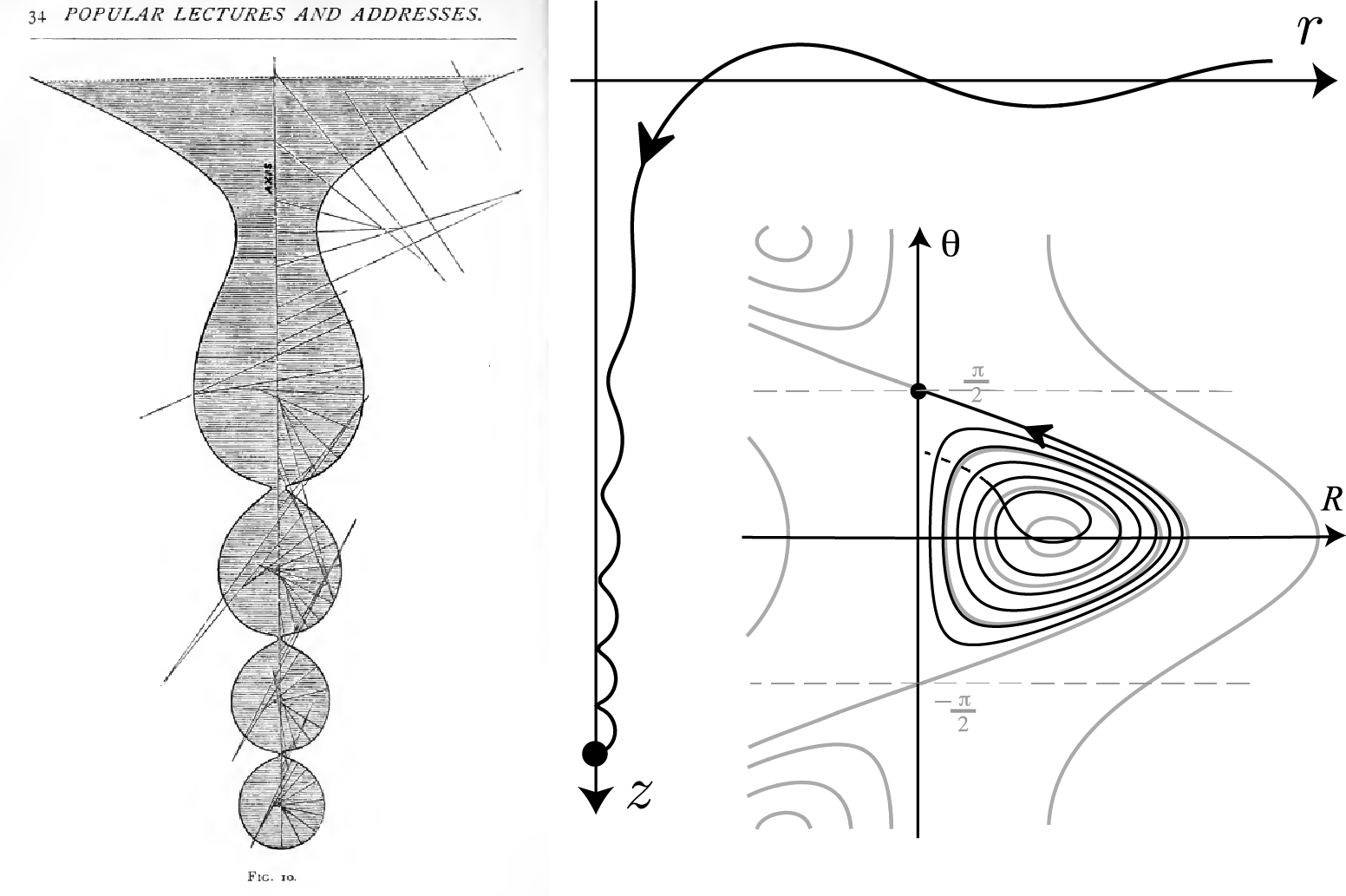}
    \caption{Left: pendent drop with four nodes, drawn by John Perry, a collaborator of Lord Kelvin, and reproduced from \cite[p. 34]{kelvin1891}. Center and right: numerical computation (reproduced from \cite[p. 1850]{rierarisler2002}) of the profile of a (finite) pendent drop (with six nodes, in the case $n=2$) and corresponding behaviour in the $(R,\theta)$-plane. The tip of the drop corresponds to the final point $(\pi/2,0)$ of the trajectory in the $(R,\theta)$-plane.}
    \label{fig_drop_phase}
\end{figure}

Expression of $d\theta/dz$ in system~(\ref{syst_rz}) suggests the following change of variables: 
\[
(r,z) \mapsto (R,Z) = \Bigl(zr, \frac{z^2}{2}\Bigr)\,,
\]
which corresponds to a blow-up of factor $z$ at depth $z$, both in the directions of $r$ and $z$, see figures \ref{fig_def_theta} and \ref{fig_drop_phase}. Rewritten using $(R,Z)$ variables, system~(\ref{syst_rz}) becomes:
\begin{equation}
\label{syst_RZ}
\left\{
\begin{aligned}
dR/dZ &= -\tan\theta+\frac{R}{2Z} \\
d\theta/dZ &= (n-1)\Bigl(\frac{1}{\cos\theta}-\frac{1}{R}\Bigr)
\end{aligned}
\right. 
\quad\mbox{where}\quad
Z>0
\,,\quad
R>0
\,,\quad
\theta\in\Bigl(-\frac{\pi}{2},\frac{\pi}{2}\Bigr)\,,
\end{equation}
which asymptotically reduces, when $Z$ approaches $+\infty$, to the two-dimensional autonomous differential system
\begin{equation}
\label{syst_RZ_lim}
\left\{
\begin{aligned}
dR/dZ &= -\tan\theta \\
d\theta/dZ &= (n-1)\Bigl(\frac{1}{\cos\theta}-\frac{1}{R}\Bigr)
\end{aligned}
\right. 
\end{equation}
\begin{figure}[!htbp]
	\centering
    \includegraphics[width=0.5\textwidth]{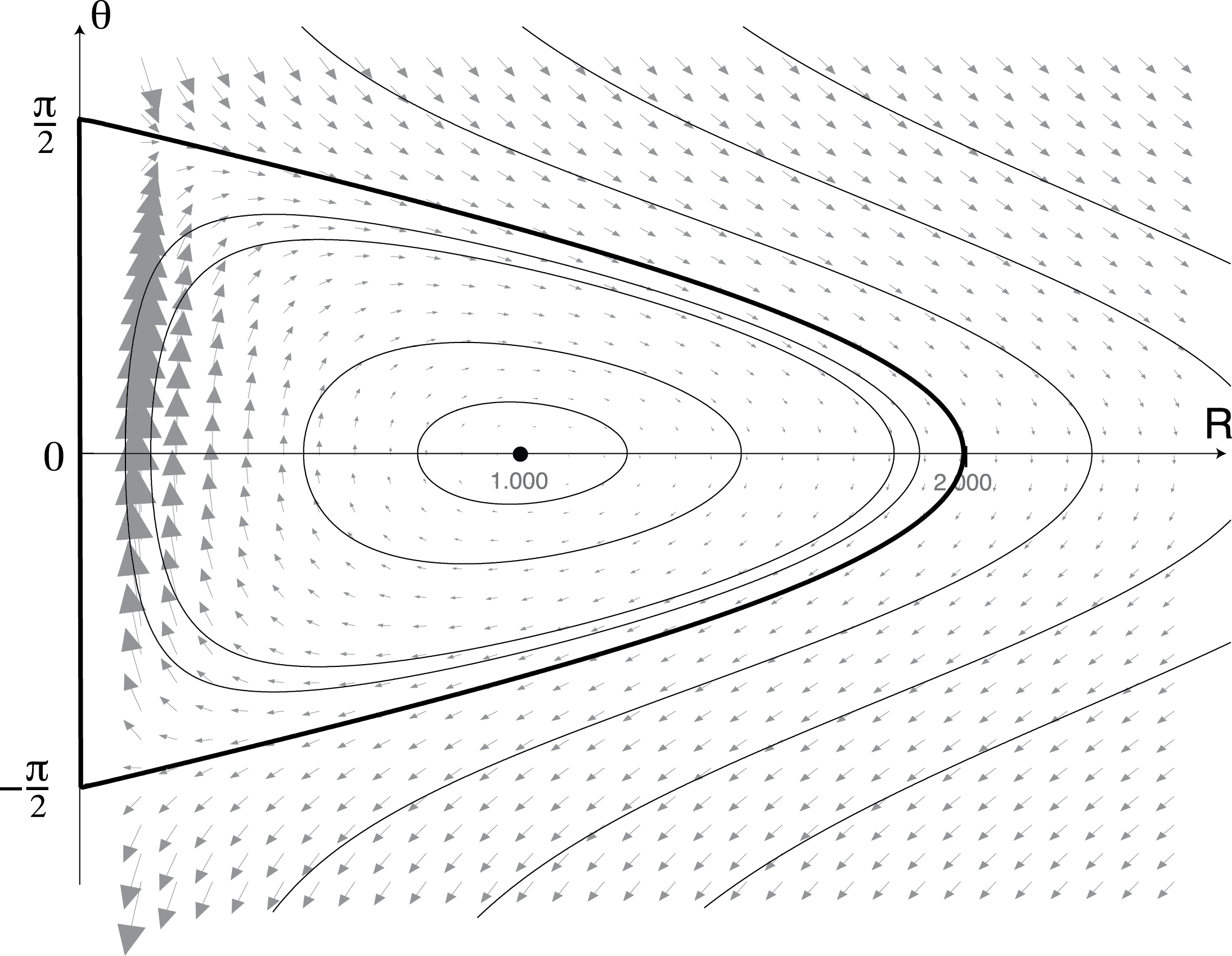}
    \caption{Phase portrait of system~(\ref{syst_RZ_lim}) (the time variable is $t$ with $dZ=-\cos\theta\,dt$, instead of $Z$) in the case $n=2$.}
    \label{fig_phase_asympt}
\end{figure}
for which the quantity
\[
\hhh(R,\theta)= R^{n-1}\Bigl(\cos\theta-\frac{n-1}{n}R\Bigr)
\]
is conserved (see figure~\ref{fig_phase_asympt}). Indeed,
\begin{equation}
\label{dH}
\left\{
\begin{aligned}
\partial_R \hhh &= (n-1) R^{n-2}(\cos\theta -R)\\
\partial_\theta \hhh &= - R^{n-1}\sin\theta
\end{aligned}
\right. 
\qquad\mbox{thus}\qquad
\left\{
\begin{aligned}
dR/dZ &= \frac{1}{R^{n-1}\cos\theta} \partial_\theta \hhh \\
d\theta/dZ &= - \frac{1}{R^{n-1}\cos\theta} \partial_R \hhh
\end{aligned}
\right. 
\end{equation}
which shows that, up to a scale change in $Z$, system~(\ref{syst_RZ_lim}) is nothing but the Hamiltonian system deriving from $\hhh(.,.)$ by the usual symplectic form on $\rr^2$. 

The phase portrait of system~(\ref{syst_RZ_lim}) is shown (in the case $n$ equals $2$) on figure~\ref{fig_phase_asympt}. It admits an elliptic equilibrium at $(1,0)$ and no other equilibrium in $(0,+\infty)\times(-\pi/2,\pi/2)$. The Hessian matrix of $\hhh$ at $(1,0)$ reads:
\[
\begin{pmatrix}
-(n-1) & 0 \\ 0 & -1
\end{pmatrix}
\]
This matrix being negative definite, the point $(1,0)$ is a strict local maximum point of $\hhh$. You could actually check that this point is also the strict global maximum point of $\hhh$ on $[0,+\infty)\times(-\pi,\pi)$. 

As was extensively studied by Riera and the author in \cite{rierarisler2002}, the effect of the additional ``perturbation'' term $R/(2Z)$ in the expression of $dR/dZ$ in ``full'' system~(\ref{syst_RZ}) is to take the solution away from the maximum point $(1,0)$ of $\hhh(.,.)$ when $Z$ increases (in other words, to decrease the value of $\hhh(R,\theta)$ along the solution). See figure~\ref{fig_drop_phase}. 
Using this repulsive effect of the perturbation term, the following existence result was proved in \cite{rierarisler2002} (once repulsion is quantified, it follows from an elementary argument, namely Cantor's intersection theorem). This result is almost identical to the existence result established by Concus and Finn in \cite{concusfinn1975a,concusfinn1975b}.
\begin{theorem}
\label{exist_thm}
There exists (at least) one solution $Z\mapsto \bigl(R(Z),\theta(Z)\bigr)$ of system~(\ref{syst_RZ}) that is defined on $(0,+\infty)$, and every such solution satisfies:
\[
\bigl(R(Z),\theta(Z)\bigr)\rightarrow(1,0)
\quad\mbox{when}\quad
Z\rightarrow+\infty\,.
\]
\end{theorem}
\section{Sketch of the proof}
\label{sec_sketch}
As can be expected as this stage, the proof of uniqueness of the singular solution relies on the same ``repulsion'' effect displayed on figure~\ref{fig_drop_phase}. We shall assume that two distinct singular solutions exist, then we shall prove by an elementary perturbation argument that the resulting effect of the ``additional'' term $R/(2Z)$ in system~(\ref{syst_RZ}) is to increase the distance between those two solutions, at least for $Z$ sufficiently large, a contradiction with the fact that both must converge to $(1,0)$ when $Z$ approaches $+\infty$. 

The sole obstacle to overcome is that the perturbation argument requires an a priori estimate on singular solutions slightly more precise than: ``$\bigl(R(Z),\theta(Z)\bigr)$ approaches $(1,0)$ when $Z$ approaches $+\infty$''. What will be required is the following fist order approximation:
\begin{equation}
\label{1st_ord_app}
\left\{
\begin{aligned}
R(Z) &= 1 + o_{Z\rightarrow+\infty}\Bigl(\frac{1}{Z}\Bigr) \\
\theta(Z) &= \frac{Z}{2} + o_{Z\rightarrow+\infty}\Bigl(\frac{1}{Z}\Bigr)
\end{aligned}
\right. 
\end{equation}

The remaining of the paper is organized as follows. In section~\ref{sec_formal} we shall describe the (unique) formal singular solution (as an asymptotic expansion in powers of $1/Z$), and observe that it matches first order expansion~(\ref{1st_ord_app}). In section~\ref{sec_perturbation} we shall carry on the main perturbation argument and, assuming that~(\ref{1st_ord_app}) holds, we shall prove uniqueness of the singular solution. In section~\ref{sec_matching}, the very same argument will be used to prove~(\ref{1st_ord_app}) --- by the way we shall actually prove that the singular solution matches the asymptotic expansion of section~\ref{sec_formal} up to every order.
\section{Formal singular solution (asymptotic expansion)}
\label{sec_formal}
In the following proposition, a \emph{formal solution} means an expansion of the form~(\ref{expr_form_sol}) below such that, when it is injected in system~(\ref{syst_RZ}) where the terms $\tan\theta$ and $\cos\theta$ are replaced by their expansions in power series of $\theta$, then the two equations of~(\ref{syst_RZ}) hold at every order in $Z^{-1}$. 
\begin{proposition}  
\label{prop_form_sol}
System~(\ref{syst_RZ}) admits a unique formal solution of the form   
\begin{equation}
\label{expr_form_sol}
\left\{
\begin{aligned}
R(Z) &= \sum_{k\in\nn}\frac{R_k}{Z^k} \\
\theta(Z) &= \sum_{k\in\nn}\frac{\theta_k}{Z^k}
\end{aligned}
\right. 
\quad\mbox{with}\quad
(R_k,\theta_k)\in\rr^2
\quad\mbox{and}\quad
R_0>0
\quad\mbox{and}\quad
-\frac{\pi}{2}<\theta_0<\frac{\pi}{2}\,.
\end{equation}
Moreover, 
\begin{itemize}
\item expansion $\sum_{k\in\nn}R_k / Z^k$ is even with respect to $Z$ (that is $R_k=0$ if $k$ is odd);
\item expansion $\sum_{k\in\nn}\theta_k / Z^k$ is odd with respect to $Z$ (that is $\theta_k=0$ if $k$ is even);
\end{itemize}
\end{proposition}
\proof
Replacing in system~(\ref{syst_RZ}) the quantities $R(Z)$ and $\theta(Z)$ by their formal expansions~(\ref{expr_form_sol}) gives:
\begin{align}
\sum_{k\ge2} -\frac{(k-1)R_{k-1}}{Z^k} &= -\tan\biggl(\sum_{k\ge0}\frac{\theta_k}{Z^k}\biggr) + \sum_{k\ge1}\frac{R_{k-1}}{2Z^k}\,, \label{form_dRdZ}\\
\sum_{k\ge2} -\frac{(k-1)\theta_{k-1}}{Z^k} &= (n-1)\biggl(\frac{1}{\cos\bigl(\sum_{k\ge0}\theta_k/Z^k\bigr)} - \frac{1}{\sum_{k\ge0}R_k/Z^k}\biggr)\,.\label{form_dthetadZ}
\end{align}
Thus,
\begin{itemize}
\item at order $k=0$, equation~(\ref{form_dRdZ}) yields $\theta_0 = 0$ and equation~(\ref{form_dthetadZ}) yields $R_0 = 1$;
\item at order $k=1$, equation~(\ref{form_dRdZ}) yields $\theta_1 = 1/2$ and equation~(\ref{form_dthetadZ}) yields $R_0 = 0$.
\end{itemize}
More generally, at every order $k$ in $\nn$, 
\begin{itemize}
\item equation~(\ref{form_dRdZ}) provides an expression of $\theta_k$ depending on $R_0$, …, $R_{k-1}$ and $\theta_0$, …, $\theta_{k-1}$;
\item equation~(\ref{form_dthetadZ}) provides an expression of $R_k$ depending on $R_0$, …, $R_{k-1}$ and $\theta_0$, …, $\theta_{k-1}$.
\end{itemize}
This proves the uniqueness of the expansion~(\ref{expr_form_sol}) satisfying system~(\ref{syst_RZ}). 
 
Now assume that this expansion of $R(Z)$ is even up to order $k-1$ and this expansion of $\theta(Z)$ is odd up to some order $k-1$. Then, 
\begin{itemize}
\item if $k$ is even then the expression of $\theta_k$ provided by equation~(\ref{form_dRdZ}) equals zero; this follows from the facts that $R_{k-1}=0$ and that tangent function is odd;
\item if $k$ is odd then the expression of $R_k$ provided by equation~(\ref{form_dthetadZ}) equals zero; this follows from the facts that $\theta_{k-1}=0$ and that cosine function is even.
\end{itemize}
Thus expansion of $R(Z)$ ($\theta(Z)$) remains even (respectively, odd) up to order $k$. This completes the proof. 
\qed
\section{Perturbation argument}
\label{sec_perturbation}
Let 
\[
Z\mapsto \bigl(R(Z),\theta(Z)\bigr)
\quad\mbox{and}\quad
Z\mapsto \bigl(\tilde R(Z),\tilde \theta(Z)\bigr) 
\]
denote two singular solutions of system~(\ref{syst_RZ}) (that is, defined up to $+\infty$ in $Z$); according to theorem~\ref{exist_thm}, both must approach $(1,0)$ when $Z$ approaches $+\infty$). Let us write:
\[
\begin{aligned}
\rho(Z) &= \tilde R(Z) - R(Z)\,, \\
\varphi(Z) &= \tilde \theta(Z) - \theta(Z)\,.
\end{aligned}
\]
Since the two solutions share the same limit when $Z$ approaches $+\infty$, it follows that:
\[
\bigl(\rho(Z),\varphi(Z)\bigr)\rightarrow (0,0)
\quad\mbox{when}\quad
Z\rightarrow+\infty\,.
\]
To enforce the contradiction argument sketched in section \ref{sec_sketch}, we would like to show that, if these two solutions differ, then the distance between them --- that is, the size of the pair $\bigl(\rho(Z), \varphi(Z)\bigr)$ --- increases with $Z$. Since $(1,0)$ is a strict local maximum point of $\hhh$, a possible approach for that is to show that the function
\[
h(Z) = \hhh\bigl(1+\rho(Z), \varphi(Z)\bigr)
\]
decreases with $Z$. Basically,
\begin{equation}
\label{dhdZ_step_0}
\frac{dh}{dZ} =\partial_R \hhh \bigl(1+\rho(Z), \varphi(Z)\bigr) \cdot\, \frac{d\rho}{dZ}+\partial_\theta \hhh \bigl(1+\rho(Z), \varphi(Z)\bigr) \cdot\, \frac{d\varphi}{dZ}
\,.
\end{equation}
From expressions~(\ref{dH}) of $\partial_R \hhh$ and $\partial_\theta \hhh$, it follows that:
\begin{align}
\partial_R \hhh\bigl(1+\rho(Z), \varphi(Z)\bigr) & = (n-1) (1+\rho)^{n-2} (\cos\varphi -1 -\rho) \nonumber\\
& = - (n-1)\Bigl( \rho + o_{Z\rightarrow+\infty}\bigl(\abs{(\rho,\varphi)}\bigr)\Bigr)\,, \label{part_R_H} \\
\partial_\theta \hhh\bigl(1+\rho(Z), \varphi(Z)\bigr) & = - (1+\rho)^{n-1} \sin\varphi \nonumber\\
& = - \Bigl( \varphi + o_{Z\rightarrow+\infty}\bigl(\abs{(\rho,\varphi)}\bigr)\Bigr)\label{part_th_H}
\end{align}
with the notation
\[
\abs{(\rho,\varphi)} = \sqrt{\rho^2+\varphi^2}\,.
\]
Let us define functions $\varepsilon_{\rho}(Z)$ and $\varepsilon_{\varphi}(Z)$ by:
\[
\left\{
\begin{aligned}
d\rho/dZ &= -\tan\varphi+ \varepsilon_{\rho} \\
d\varphi/dZ &= (n-1)\Bigl(\frac{1}{\cos\varphi}-\frac{1}{1+\rho}\Bigr)+\varepsilon_{\varphi}
\end{aligned}
\right. 
\]
The quantities $\varepsilon_{\rho}(Z)$ and $\varepsilon_{\varphi}(Z)$ thus measure the discards in $d\rho/dZ$ and $d\varphi/dZ$ with respect to what these two derivatives would be if the pair $\bigl(1+\rho(Z),\varphi(Z)\bigr)$ was (exactly) governed by the asymptotic differential system~(\ref{syst_RZ_lim}). Since $\hhh(.,.)$ is a conserved quantity for this system, expression~(\ref{dhdZ_step_0}) of $dh/dZ$ then reduces to
\[
\frac{dh}{dZ} =\partial_R \hhh \cdot\, \varepsilon_{\rho}+\partial_\theta \hhh \cdot\, \varepsilon_{\varphi}
\]
and, substituting expansions~(\ref{part_R_H}) and~(\ref{part_th_H}) of $\partial_R \hhh $ and $\partial_\theta \hhh$ in this expression, it follows that
\begin{equation}
\label{dhdZ_step_1}
\frac{dh}{dZ} = - (n-1)\Bigl( \rho + o_{Z\rightarrow+\infty}\bigl(\abs{(\rho,\varphi)}\bigr)\Bigr)  \cdot \varepsilon_{\rho} - \Bigl( \varphi + o_{Z\rightarrow+\infty}\bigl(\abs{(\rho,\varphi)}\bigr)\Bigr) \cdot \varepsilon_{\varphi} ,.
\end{equation}
This leads us to compute the quantities $\varepsilon_{\rho}(Z)$ and $\varepsilon_{\varphi}(Z)$. It follows from system~(\ref{syst_RZ}) that
\[
\begin{aligned}
\varepsilon_{\rho} = {} & -(\tan\tilde\theta - \tan\theta) + \frac{\rho}{2Z} + \tan\varphi \\
= {} & \frac{\rho}{2Z} -\tan\varphi \tan\theta \tan(\theta+\varphi)
\,, \\
\varepsilon_{\varphi} = {} & (n-1) \biggl( \Bigl(\frac{1}{\cos(\theta+\varphi)}-\frac{1}{\cos\theta}\Bigr)-\Bigl(\frac{1}{\cos\varphi}-1\Bigr)
+
\Bigl(-\frac{1}{R+\rho} +\frac{1}{R}\Bigr) + \Bigl(\frac{1}{1+\rho}-1\Bigr) \biggr) \\
= {} &  (n-1)\biggl( \frac{\sin\varphi\sin\theta}{\cos\theta\cos(\theta+\varphi)}  \\
& + \frac{1-\cos\varphi}{\cos(\theta+\varphi)}\Bigl((1-\cos\theta) +\frac{\sin\varphi\sin\theta}{\cos\varphi}\Bigr) + \rho(1-R)\frac{1+R+\rho}{R(R+\rho)(1+\rho)}\biggr)\,.
\end{aligned}
\]
At first glance these expressions look a bit intricate, but if we make the additional hypothesis~(\ref{1st_ord_app}), that is if we assume that
\[
\begin{split}
R(Z) &= 1 + o_{Z\rightarrow+\infty}\Bigl(\frac{1}{Z}\Bigr)\,, \\
\theta(Z) &= \frac{1}{2Z} + o_{Z\rightarrow+\infty}\Bigl(\frac{1}{Z}\Bigr)\,,
\end{split}
\]
then expressions of $\varepsilon_{\rho}$ and $\varepsilon_{\varphi}$ above yield:
\[
\begin{aligned}
\varepsilon_{\rho}(Z) &= \frac{\rho}{2Z} + \frac{o_{Z\rightarrow+\infty}\bigl(\abs{(\rho,\varphi)}\bigr)}{Z}\,,  \\
\varepsilon_{\varphi}(Z) &= (n-1) \frac{\varphi}{2Z} + \frac{o_{Z\rightarrow+\infty}\bigl(\abs{(\rho,\varphi)}\bigr)}{Z}\,.
\end{aligned}
\]
Expansion~(\ref{dhdZ_step_1}) of $dh/dZ$ thus becomes
\[
\frac{dh}{dZ} = -(n-1)\frac{\rho^2+\varphi^2}{2Z}\bigl(1 + o_{Z\rightarrow+\infty}(1)\bigr)\,,
\]
and as a consequence, for all $Z$ sufficiently large, 
\begin{equation}
\label{dhdZ_final}
\frac{dh}{dZ}\le 0.
\end{equation}

On the other hand, since $\bigl(\rho(Z),\varphi(Z)\bigr)$ approaches $(0,0)$ when $Z$ approaches $+\infty$, it follows that:
\[
h(Z)\rightarrow \hhh(1,0) 
\quad\mbox{when}\quad
Z\rightarrow+\infty\,.
\]
Since $(1,0)$ is a strict local maximum point of $\hhh$, this shows, together with inequality~(\ref{dhdZ_final}), that, for all $Z$ sufficiently large, $h(Z)=\hhh(1,0)$, in other words:
\[
\bigl(\rho(Z),\varphi(Z)\bigr)=(0,0)\,. 
\]
The two singular solutions $(R,\theta)$ and $(\tilde{R},\tilde{\theta})$ are thus equal. This proves the desired uniqueness under hypothesis~(\ref{1st_ord_app}).

To complete the proof of theorem~\ref{main_theorem}, the sole remaining thing to prove is that hypothesis~(\ref{1st_ord_app}) holds for every singular solution. 
\section{Matching with asymptotic expansion}
\label{sec_matching}
The following proposition states that every singular solution coincides, up to every order, with the formal singular solution defined in proposition~\ref{prop_form_sol}. 
\begin{proposition}    
\label{prop_approx} 
For every solution $Z\mapsto \bigl(\tilde R(Z),\tilde \theta(Z)\bigr)$ of system~(\ref{syst_RZ}), defined on $(0,+\infty)$, and for all $k$ in $\nn$, the following expansions hold:
\[
\begin{aligned}
\tilde R(Z) &= \sum_{j=0}^k \frac{R_j}{Z^j} + o_{Z\rightarrow+\infty}\Bigl(\frac{1}{Z^k}\Bigr)\,, \\
\tilde\theta(Z) &= \sum_{j=0}^k \frac{\theta_j}{Z^j} + o_{Z\rightarrow+\infty}\Bigl(\frac{1}{Z^k}\Bigr)\,.
\end{aligned}
\]
\end{proposition}
\proof
Let $Z\mapsto \bigl(\tilde R(Z),\tilde\theta(Z)\bigr)$ denote a solution of of system~(\ref{syst_RZ}) defined on $(0,+\infty)$. According to theorem~\ref{exist_thm}, the pair $\bigl(\tilde R(Z),\tilde\theta(Z)\bigr)$ approaches $(1,0)$ when $Z$ approaches $+\infty$. 
Since $R_0=1$ and $\theta_0=0$, this proves the desired result for $k=0$. 

Now take $k$ in $\nn^*$ and let us consider the following two functions of $Z$: 
\[
R(Z) = \sum_{j=0}^k \frac{R_j}{Z^j}
\quad\mbox{and}\quad
\theta(Z) = \sum_{j=0}^k \frac{\theta_j}{Z^j}
\]
(note that these two functions depend on the choice of $k$).
Then, the definitions of the coefficients $R_j$ and $\theta_j$ in the proof of proposition~\ref{prop_form_sol} show that the following estimates hold: 
\begin{equation}
\label{dapprox_dZ}
\left\{
\begin{aligned}
d R/dZ &= -\tan\theta+\frac{R}{2Z} +\ooo_{Z\rightarrow+\infty}\Bigl(\frac{1}{Z^{k+1}}\Bigr) \\
d \theta/dZ &= (n-1)\Bigl(\frac{1}{\cos\theta}-\frac{1}{R}\Bigr)+\ooo_{Z\rightarrow+\infty}\Bigl(\frac{1}{Z^{k+1}}\Bigr)
\end{aligned}
\right. 
\end{equation}
Moreover, since $k\ge1$, hypothesis~(\ref{1st_ord_app}) holds for $(R,\theta)$:
\[
\begin{split}
R(Z) &= 1+o_{Z\rightarrow+\infty}\Bigl(\frac{1}{Z}\Bigr)\,, \\
\theta(Z) &= \frac{1}{2Z}+o_{Z\rightarrow+\infty}\Bigl(\frac{1}{Z}\Bigr)\,. 
\end{split}
\]

To complete the proof, we will apply exactly the same approach (and adopt the same notation) as in previous section, the only difference being the remaining terms $\ooo(1/Z^{k+1})$ on the right-hand side of system~(\ref{dapprox_dZ}) above. Let us write:
\[
\begin{aligned}
\rho (Z)&= \tilde R(Z) - R(Z)\,, \\
\varphi(Z) &= \tilde \theta(Z) - \theta(Z)\,, \\
h(Z) &= \hhh\bigl(1+\rho(Z), \varphi(Z)\bigr)\,.
\end{aligned}
\]
Expansions~(\ref{part_R_H}) and~(\ref{part_th_H}) of $\partial_R\hhh$ and $\partial_\theta\hhh$ hold unchanged. 
Let use define functions $\varepsilon_{\rho}(Z)$ and $\varepsilon_{\varphi}(Z)$ by:
\[
\left\{
\begin{aligned}
d\rho/dZ &= -\tan\varphi+ \varepsilon_{\rho} \\
d\varphi/dZ &= (n-1)\Bigl(\frac{1}{\cos\varphi}-\frac{1}{1+\rho}\Bigr)+\varepsilon_{\varphi}
\end{aligned}
\right. 
\]
Expressions of $\varepsilon_{\rho}$ and $\varepsilon_{\varphi}$ are identical to those of previous section, except the fact that they comprise additional terms $\ooo(1/Z^{k+1})$. Thus,
\[
\begin{aligned}
\varepsilon_{\rho} &= \frac{\rho}{2Z} + \frac{o_{Z\rightarrow+\infty}\bigl(\abs{(\rho,\varphi)}\bigr)}{Z} + \ooo_{Z\rightarrow+\infty}\Bigl(\frac{1}{Z^{k+1}}\Bigr)\,,  \\
\varepsilon_{\varphi} &= (n-1) \frac{\varphi}{2Z} + \frac{o_{Z\rightarrow+\infty}\bigl(\abs{(\rho,\varphi)}\bigr)}{Z} + \ooo_{Z\rightarrow+\infty}\Bigl(\frac{1}{Z^{k+1}}\Bigr)\,,
\end{aligned}
\]
and it follows from unchanged expansion~(\ref{dhdZ_step_1}) that:
\begin{align}
\frac{dh}{dZ} &= -(n-1)\frac{\rho^2+\varphi^2}{2Z}\bigl(1 + o_{Z\rightarrow+\infty}(1)\bigr) + \ooo_{Z\rightarrow+\infty}\Bigl(\frac{\abs{(\rho,\varphi)}}{Z^{k+1}}\Bigr) \nonumber\\
& = -(n-1)\frac{\abs{(\rho,\varphi)}}{2Z}\biggl(\abs{(\rho,\varphi)}\bigl(1+o_{Z\rightarrow+\infty}(1)\bigr) + \ooo_{Z\rightarrow+\infty}\Bigl(\frac{1}{Z^{k}}\Bigr) \biggr)\,. \label{dhdZ_perturb}
\end{align}

We are going to deduce from this last estimate that
\[
\abs{(\rho,\varphi)} = \ooo_{Z\rightarrow+\infty} \Bigl(\frac{1}{Z^{k}}\Bigr)\,.
\]
Indeed, according to expression~(\ref{dhdZ_perturb}) of $dh/dZ$, there exist positive quantities $C_1$ and $Z_0$ such that, for all $Z$ superior or equal to $Z_0$,
\[
\abs*{\bigl(\rho(Z),\varphi(Z)\bigr)} \ge \frac{C_1}{Z^k} \Longrightarrow \frac{dh}{dZ} \le 0. 
\]
Besides, since $(1,0)$ is a strict local maximum point of $\hhh$, since the Hessian matrix of $\hhh$ at $(0,1)$ is negative definite, and since $\hhh(1,0)$ equals $1/n$, the quantities $\abs{(\rho,\varphi)}$ and $\sqrt{1/n - h(Z)}$ do not differ, if $Z$ is large, by more than a fixed multiplicative factor. As a consequence, there exists a positive quantity $C_2$ such that, for all $Z$ greater than or equal to $Z_0$,
\[
\sqrt{\frac{1}{n} - h(Z)} \ge \frac{C_2}{Z^k} \Longrightarrow \frac{dh}{dZ}(Z) \le 0\,,
\]
in other words:
\begin{equation}
\label{implic_h}
h(Z)\le \frac{1}{n} - \Bigl(\frac{C_2}{Z^{k}}\Bigr)^2 \Longrightarrow \frac{dh}{dZ}(Z) \le 0\,.
\end{equation}
But since $h(Z)$ approaches $1/n$ when $Z$ approaches $+\infty$, the left-hand inequality of implication~(\ref{implic_h}) can actually not occur if $Z\ge Z_0$. Indeed, by contradiction, if there existed a quantity $Z_1$ superior or equal to $Z_0$ such that
\[
h(Z_1)\le \frac{1}{n} - \Bigl(\frac{C_2}{Z1^{k}}\Bigr)^2\,,
\]
then according to implication~(\ref{implic_h}), it would follow that, for all $Z$ superior or equal to $Z_1$,
\[
h(Z)\le \frac{1}{n} - \Bigl(\frac{C_2}{Z_1^{k}}\Bigr)^2 \le \frac{1}{n} - \Bigl(\frac{C_2}{Z^{k}}\Bigr)^2\,,
\]
a contradiction with the fact that $h(Z)$ approaches $1/n$ when $Z$ approaches $+\infty$. 

We thus proved that, for all $Z$ sufficiently large, 
\[
h(Z)\ge 1/n - C_2^2/Z_1^{2k}
\qquad\mbox{or equivalently}\qquad
\sqrt{1/n - h(Z)}\le C_2/Z^k\,, 
\]
and this finally yields the desired estimate:
\[
\abs{(\rho,\varphi)} = \ooo_{Z\rightarrow+\infty}\Bigl(\frac{1}{Z^{k}}\Bigr)\,.
\]
Since the positive integer $k$ was arbitrary, this finishes the proof of proposition~\ref{prop_approx}.
\qed

In view of the comments in section~\ref{sec_perturbation}, proof of theorem~\ref{main_theorem} is complete. 
\subsubsection*{Acknowledgements}
I thank Pierre Coullet who introduced me to the capillary equation and this perturbation approach, and Christophe Riera who carried out the numerical simulations and with whom I had many fruitful discussions. 
\printbibliography 
%
%
\end{document}